\documentclass{article}[12pt]

\usepackage{latexsym}
\usepackage{amssymb}
\usepackage{amsmath,amsfonts}
\usepackage{graphicx}
\usepackage{tikz}
\usetikzlibrary{arrows}

\newtheorem{theorem}{Theorem}[section]
\newtheorem{prop}[theorem]{Proposition}
\newtheorem{lemma}{Lemma}[section]
\newtheorem{corollary}{Corollary}[section]
\newtheorem{example}{Example}[section]
\newtheorem{definition}{Definition}[section]

\newtheorem{conjecture}[theorem]{Conjecture}

\newtheorem{question}[theorem]{Question}

\def \dd{\hfill \baseb \vskip .5cm}
\def \d{{\noindent \it Proof. } }

\def \zk{{\mathbb Z}_k}
\def \N{{\mathbb N}}
\def \qq{{\mathbb Q}}
\def\A{{\mathcal A}}

\def\ex{\begin{example}}
\def\eex{\end{example}}
\def\exx{\end{example}}
\def\t{\begin{theorem}}
\def\tt{\end{theorem}}
\def\D{\begin{definition}}
\def\DD{\end{definition}}
\def\l{\begin{lemma}}
\def\ll{\end{lemma}}
\def\c{\begin{corollary}}
\def\cc{\end{corollary}}
\def\cj{\begin{conjecture}}
\def\cjj{\end{conjecture}}
\def\e{\begin{equation}}
\def\ee{\end{equation}}
\def\p{\begin{prop}}
\def\pp{\end{prop}}
\def\qu{\begin{question}}
\def\qqu{\end{question}}

\def \Q{{\mathcal Q}}

\def \Dn{{\mathcal D}_n}

\def \Gn{{\mathcal G}_n}

\def \Tn{{\mathcal T}_n}

\newcommand{\baseb}{\hfill \rule{2mm}{2mm}}

\begin{document}

\baselineskip 15pt

\title{{\bf Cordiality of Digraphs}\thanks{AMS Classification number: 05C20, 05C38, 05C78
} \thanks{Key words and phrases: tournament, wheel graph, fan graph, $(2,3)$-cordial)
}}
\author{LeRoy B. Beasley, Manuel A. Santana,   Jonathan Mousley\\ and David E. Brown\\ \\ Department of Mathematics and Statistics \\Utah State University \\ Logan, Utah 84322-3900, U.S.A}
%\date{}

\maketitle

%\pagebreak
\begin{abstract}
A $(0,1)$-labelling of a set is said to be {\em friendly} if approximately one half the elements of the set are labelled 0 and one half labelled 1.  Let $g$ be a labelling of the edge set of a graph that is induced by a labelling $f$ of the vertex set.  If both $g$ and $f$ are friendly then $g$ is said to be a {\em cordial} labelling of the graph.  We extend this concept to directed graphs and investigate the cordiality of sets of directed graphs.  We investigate a specific type of cordiality on digraphs, a restriction  of quasigroup-cordiality called  $(2,3)$-cordiality.   A directed graph is $(2,3)$-cordial if there is a  friendly labelling $f$ of the vertex set which induces a $(1,-1,0)$-labelling of the arc set $g$ such  that  about one third of the arcs are labelled 1, about one third labelled -1 and about one third labelled 0.   In particular we determine which tournaments are $(2,3)$-cordial,  which orientations of the $n$-wheel are $(2,3)$-cordial, and which orientations of the $n -$fan are $(2,3)$-cordial.
\end{abstract}

\section{Introduction}

The study of cordial graphs began in 1987 with an article by I. Cahit \cite{C}: "Cordial Graphs: A Weaker Version of Graceful and Harmonious Graphs".  In 1991, Hovey \cite{H} generalized this concept to $\A$-cordial graphs where $\A$ is an abelian group.  A further generalization,  one that included cordiality of directed graphs, appeared  in 2012 with an article by Pechenik and  Wise  \cite{PW}, where the $\A$ was allowed to be any quasi group, not necessarily Abelian.  We modify this concept to $(A,\A)$-cordial digraphs where $A$ is a subset of the quasigroup $\A$. 

Let $G=(V,E)$ be an undirected graph with vertex set $V$ and edge set $E$.  A $(0,1)$-labelling of the vertex set is a mapping $f:V\to \{0,1\}$ and is said to be {\em friendly} if approximately one half of the vertices are labelled 0 and the others labelled 1.  An induced labelling of the edge set  is a mapping $g:E\to \{0,1\}$ where for an edge $uv, g(uv)= \hat{g}(f(u),f(v)$ for some $\hat{g}:\{0,1\}\times\{0,1\}\to \{0,1\}$ and is said to be cordial if $f$ is friendly and about one half the edges of $G$ are labelled 0 and the others labelled 1.  A graph, $G$,  is called {\em cordial} if there exists a cordial induced labelling of the edge set of $G$.  In this article, as in \cite{B}, we define a cordial labelling of directed graphs that is not  merely a cordial labelling of the underlying undirected graph.  

A specific type of $(A,\A)$-cordial digraph is a $(2,3)$-cordial digraph defined by Beasley in \cite{B}.  Let $D=(V,A)$ be a directed graph with vertex set $V$ and arc set $A$.  Let $f:V\to \{0,1\}$ be a friendly vertex labelling and let $g$ be the  induced labelling of the arc set,  $g:A\to \{0,1,-1\}$ where for an arc $\overrightarrow{uv}, g(\overrightarrow{uv})=f(v)-f(u)$.  The labellings $f$ and $g$  are  {\em $(2,3)$-cordial} if $f$ is friendly and about one third the arcs of $D$ are labelled 1, one third are labelled -1 and one third labelled 0.   A digraph, $D$,  is called {\em $(2,3)$-cordial} if there exists  $(2,3)$-cordial  labellings $f$ of the vertex set and $g$  of the arc set of $D$. 

Note that here and what follows, the term ``about'' when talking about fractions of a quantity we shall mean as close is possible in integral arithmetic, so about half of 9 is either 4 or 5, but not 3 or 6.

\section{Preliminaries}

\D  A {\em quasigroup} is a set $\Q$  with binary operation $\circ $ such that
 given any $a,b\in\Q$ there are $x,y\in \Q$ such that $a\circ x=b$ and $y\circ a = b$. \DD

\noindent {\bf \Large Fact:}  {\em All two element quasigroups are Abelian.}

\d Suppose that $\Q=\{a,b\}$ is a quasi group with binary operation $\circ$.  Then, there are  $x,y\in \Q$ such that $a\circ x=a$ and $y\circ a = a$. If $x=y=b$ then $\Q$ is Abelian.  Otherwise, we must have $a\circ a = a$.  Similarly either $\Q$ is abelian or $b\circ b=b$. 

Now, suppose that $a\circ a = a$ and  $b\circ b=b$.  Then  there are  $c,d\in \Q$ such that $a\circ d=b$ and $c\dot a = b$.  Since  $a\circ a = a$., we must have that both $c=b$ and $d=b$.  That is $\Q$ is Abelian.\dd 

We now formalize the terms mentioned in the introduction.  We let $\zk$ denote the set of integers $\{0,1,\dots,k\}$ with arithmetic is modulo $k$ as needed.  Further let $\zk^-$ denote the set  $\zk$ with binary operation ``-'', subtraction modulo $k$.  Clearly $\zk^-$ is a nonabelian quasigroup.

\D  A $\zk$-labelling (or simply a $k$-labelling) of a finite set, ${\cal X}$,  is a mapping $f:{\cal X}\to \zk$ and is  said to be \underline{\em friendly} if the labelling is evenly distributed over $\zk$,  that is, given any $i,j\in\zk$,  $-1\leq |f^{-1}(i)|-|f^{-1}(j)| \leq 1$ where $|{\cal X}|$ denotes the cardinality of the set ${\cal X}$.  \DD

\D Let $G=(V,E)$ be an undirected graph with vertex set $V$ and edge set $E$, and let $f$ be a friendly $(0,1)$-labelling of the vertex set $V$.   Given this friendly vertex labelling $f$, an induced $(0,1)$-labelling of the edge set  is a mapping $g:E\to \{0,1\}$ where for an edge $uv,\, g(uv)= \hat{g}(f(u),f(v))$ for some $\hat{g}:\{0,1\}\times\{0,1\}\to \{0,1\}$ and is said to be \underline{\em cordial} if $g$ is also friendly, that is about one half the edges of $G$ are labelled 0  and the others are labelled 1, or $-1\leq |g^{-1}(0)|-|g^{-1}(1)| \leq 1$.  A graph, $G$,  is called \underline{\em cordial} if there exists a induced cordial labelling of the edge set of $G$.  \DD The induced labelling $g$ in a cordial graph  is usually $g(u,v)=\hat{g}(f(u),f(v))=|f(v)-f(u)|$ \cite{C},  $g(u,v)=\hat{g}(f(u),f(v))=f(v)+f(u)$ (in ${\mathbb Z}_2$) \cite{H},  or $g(u,v)=\hat{g}(f(u),f(v))=f(v)f(u)$ (product cordiality) \cite{S}.  

In \cite{H}, Hovey introduced ${\cal A}$-friendly labellings where $\A$ is an Abelian group.  A labelling $f:V(G)\to {\cal A}$ is said to be \underline{\em ${\cal A}$-friendly } if given any $a,b\in{\cal A}$, $-1\leq |f^{-1}(b)|-|f^{-1}(a)|\leq 1$.  If $g$ is an induced edge labelling and $f$ and $g$ are both ${\cal A}$-friendly Then $g$ is said to be an ${\cal A}$-cordial labelling and $G$ is said to be ${\cal A}$-cordial.  When ${\cal A}= \zk$ we say that $G$ is $k$-cordial.  We shall use this concept with digraphs.

Given an undirected  graph or a digraph, $G$, let $\hat{G}$ denote the subgraph (or subdigraph) of $G$ induced by its nonisolated vertices.  So $\hat{G}$ never has an isolated vertex.  The need for this will become apparent in Example \ref{small}.  

In this article, we will be concerned mainly with digraphs.  We let $\Dn$ denote the set of all simple directed graphs on the vertex set $V = \{v_1,v_2,\dots,v_n\}$.  Note that the arc set of members of $\Dn$ may contain digons, a pair of arcs between two vertices each directed opposite from the other.   We shall let $\Tn$ denote the set of all subdigraphs of a tournament digraph.  So the members of $\Tn$ contain no digons.   Let $D\in\Dn$, $D=(V,A)$ where $A$ is the arc set of $D$.  Then $D$  has no loops,  and no multiple arcs.  An arc in $D$ directed from vertex $u$ to vertex $v$ will be denoted $\overrightarrow{uv}, \overleftarrow{vu}$ or by the ordered pair $(u,v)$. We also let $\Gn$ denote the set of all simple undirected graphs on the vertex set $V = \{v_1,v_2,\dots,v_n\}$.  So all members of $\Tn$ are orientations of graphs in $\Gn$.    

In \cite{PW}, Pechenik and Wise introduced quasigoup cordiality.  When the quasigroup is nonabelian, this type of cordiality is quite suitable for studying labellings of directed graphs.  In fact, if $\Q$ is a quasigroup with any binary operation $\circ$ with the property that for any $a,b\in\Q$ $a\circ b=b\circ a$ if and only if $a=b$,  we have the best situation for directed graphs.  

Now the set $\zk$ with binary operation $\circ$ where for $a,b\in\zk$, $a\circ b=b-a \mod k$ is such a quasigroup.

In our incvestigations we make one further restriction:  we will label our vertices with only a subset of $\Q$, not necessarily the whole set $\Q$:

\D Let $\Q$ be a quasigroup with binary operation $\circ$ and let $\qq$ be a subset of $\Q$. Let $D=(V,A)$ be a directed graph with vertex set $V$ and arc set $A$.   Let $f:V\to\qq$ be a friendly $\qq$-labelling of $V$ and let $g:A\to \Q$ be an induced arc labelling where for $\overrightarrow{uv}\in A$,  $g(\overrightarrow{uv})= \hat{g}(f(u),f(v))$ for some $\hat{g}:\qq\times\qq\to \Q$.  The mapping $g$ is  said to be \underline{\em $(\qq,\Q)$-cordial} if $g$ is also friendly, that is,  given any $a,b\in\Q$, $-1\leq |g^{-1}(a)|-|g^{-1}(b)| \leq 1$.  A directed graph, $D$,  is called \underline{\em $(\qq,\Q)$-cordial} if there exists a induced $(\qq,\Q)$-cordial labelling of the arc set of $D$. \DD

We now shall restrict our attention to the smallest case of $(\qq,\Q)$-cordiality that is appropriate for directed graphs, $({\mathbb Z}_2,{\mathbb Z}_3^-)$-cordiality, that defined by Beasley in \cite{B}, $(2,3)$-cordiality.

\section{$(2,3)$-orientable Digraphs.}
Let $D=(V,A)$ be a  directed graph with vertex set $V$ and arc set $A$.   Let $f:V\to\{0,1\}$ be a friendly labelling of the vertices of $D$.  As for undirected graphs, an  induced labelling of the arc set  is a mapping $g:A\to {\cal X}$ for some set ${\cal X}$ where for an arc $(u,v)=\overrightarrow{uv}, g(u,v)= \hat{g}(f(u),f(v))$ for some $\hat{g}:\{0,1\}\times\{0,1\}\to {\cal X}$.  As we are dealing with directed graphs, it would be desirable for the induced labelling to distinguish between the label of the arc $(u,v)$  and the label of the arc $(v,u)$,  otherwise, the labelling would be an induced labelling of the underlying undirected graph.    If we let ${\cal X}=\{-1,0,1\}$ and $\hat{g}(f(u),f(v))=f(v)-f(u)$ using real arithmetic, or arithmetic in ${\mathbb Z}_3$,  we have an asymmetric labelling.  In this case, if about one third of the arcs are labelled 0, about one third of the arcs are labelled 1 and about one third of the arcs are labelled -1 we say that the labelling is $(2,3)$-cordial.   Formally:

\D  Let $D\in\Tn$, $D=(V,A)$, be a digraph without isolated vertices.  Let  $f:V\to \{0,1\}$ be a friendly labelling of the vertex set $V$ of $D$.  Let $g:A\to\{1,0,-1\}$ be an induced labelling of the arcs of $D$ such that  for any $i,j\in\{1,0,-1\}$, $-1\leq |g^{-1}(i)|-|g^{-1}(j)| \leq 1$.   Such a labelling is called a \underline{\em$(2,3)$-cordial} labelling.

A digraph $D\in\Tn$ whose subgraph $\hat{G}$  can possess a $(2,3)$-cordial labelling will be called a  \underline{\em$(2,3)$-cordial} digraph.

An undirected graph $G$ is said to be \underline{\em $(2,3)$-orientable} if there exists an orientation of $G$ that is $(2,3)$-cordial. \DD

In \cite{B} the concept of $(2,3)$-cordial digraphs was introduced and paths and cycles were investigated. In \cite{SBMB} one can find  further investigation of orientations of paths and trees as well as finding the maximum number of arcs possible in a $(2,3)$-cordial digraph.  In this article we continue this investigation, showing which tournaments, which orientations of the wheel graphs, and which orientations of the fan graphs are $(2,3)$-cordial.

\D Let $D=(V,A)$ be a digraph with vertex labelling $f:V\to \{0,1\}$ and with induced arc labelling $g:A\to\{0,1,-1\}$.  Define $\Lambda_{f,g} : \Dn\to\N^3$ by $\Lambda_{f,g}(D)=(\alpha,\beta,\gamma)$ where $\alpha=|g^{-1}(1)|, \beta=|g^{-1}(-1)|,$ and $\gamma =|g^{-1}(0)|$.  \DD

Let $D\in\Tn$ and let $D^R$ be the digraph such that every arc of $D$ is reversed, so that $\overrightarrow{uv}$ is an arc in $D^R$ if and only if $\overrightarrow{vu}$ is an arc in $D$.  Let $f$ be a $(0,1)$-labelling of the vertices of $D$ and let $g(\overrightarrow{uv}) =f(v)-f(u)$ so that $g$ is a $(1,-1,0)$-labelling of the arcs of $D$.  Let $\overline{f}$ be the complementary $(0,1)$-labelling of the vertices of $D$, so that $\overline{f}(v)=0$ if and only if $f(v)=1$.  Let $\overline{g}$ be the corresponding induced arc labelling of $D$, $\overline{g}(\overrightarrow{uv}) =\overline{f}(v)-\overline{f}(u)$.

\l \label{lab} Let $D\in\Tn$ with vertex labelling $f$ and induced arc labelling $g$.  Let $\Lambda_{f,g}(D)=(\alpha,\beta,\gamma)$.   Then \begin{enumerate}\item $\Lambda_{f,g}(D^R)=(\beta,\alpha,\gamma)$. \item $\Lambda_{\overline{f},\overline{g}}(D)=(\beta,\alpha,\gamma)$, and \item $\Lambda_{\overline{f},\overline{g}}(D^R)=\Lambda_{f,g}(D).$\end{enumerate} \ll
\d If an arc is labelled 1, -1, 0 respectively then reversing the labelling of the incident vertices gives a labelling of -1, 1, 0 respectively,  If an arc $\overrightarrow{uv}$ is labelled 1, -1, 0 respectively, then $\overrightarrow{vu}$ would be labelled -1, 1, 0 respectively. \dd
 
\ex \label{small}Now, consider a graph, {  X}$_n$ in $\Gn$   consisting of three parallel edges and  n-6 isolated vertices.  Is  {  X}$_n$  $(2,3)$-orientable?  If $n=6$, the answer is no, since any friendly labelling of the six vertices would have either  no arcs labelled 0 or two arcs labelled 0.  In either case, the orientation would never be $(2,3)$-cordial. That is  {  X}$_6$ is not $(2,3)$-orientable, however  with additional vertices like  {  X}$_7$ the graph is $(2,3)$-orientable.\exx  Thus, for our investigation here, we will us the convention that a graph, $G$, is $(2,3)$-orientable/$(2,3)$-cordial  if and only if  the subgraph of $G$ induced by its nonisolated vertices, $\hat{G}$, is $(2,3)$-orientable/$(2,3)$-cordial.  

\subsection{$(2,3)$-Orientations of a Complete Graph -- Tournaments}

It is an easy exercise to show that every 3-tournament is $(2,3)$-cordial and that two of the four non isomorphic 4-tournaments are $(2,3)$-cordial.  See Figures \ref{3tour} and \ref{4-tours}.  Note that the 4-tournaments that are not $(2,3)$-cordial may require more that a cursory glance to verify that they are not $(2,3)$-cordial.

\begin{figure}[h]
\begin{center}
\begin{tikzpicture}[scale=1.5]

\tikzset{vertex/.style = {shape=circle,draw,minimum size=2em}}
\tikzset{edge/.style = {->,> = stealth',shorten >=1pt,thick}}

\node[vertex] (v1) at  (2.5,5.6) {$0$};
\node[vertex] (v2) at  (5,3.5) {$0$};
\node[vertex] (v3) at  (4,0) {$1$};
\node[vertex] (v4) at  (1,0) {$1$};
\node[vertex] (v5) at  (0,3.5) {$0$};

\draw[edge,->, line width=1.0pt] (v1) to (v2);
\draw[edge,->, line width=1.0pt] (v2) to (v3);
\draw[edge,->, line width=1.0pt] (v3) to (v4);
\draw[edge,->, line width=1.0pt] (v4) to (v5);
\draw[edge,->, line width=1.0pt] (v2) to (v5);
\draw[edge,->, line width=1.0pt] (v1) to (v3);
\draw[edge,->, line width=1.0pt] (v5) to (v3);
\draw[edge,->, line width=1.0pt] (v4) to (v1);
\draw[edge,->, line width=1.0pt] (v4) to (v2);
\draw[edge,->, line width=1.0pt] (v1) to (v5);

\path (1 ,4.5) node     (y1) {$0$};
\path (4 ,4.5) node     (y1) {$0$};
\path (2.5 ,3.6) node     (y1) {$0$};
\path (1.6,2.9) node     (y1) {$-1$};
\path (3.4,2.9) node     (y1) {$1$};
\path (0.3,1.8) node     (y1) {$-1$};
\path (4.7,1.8) node     (y1) {$1$};
\path (3.1,1.6) node     (y1) {$-1$};
\path (1.9,1.6) node     (y1) {$1$};
\path (2.5,.1) node     (y1) {$0$};

\path (4.3,-.4) node     (y1) {$v_1$};\path (.7,-.4) node     (y1) {$v_2$};\path (-.4,3.6) node     (y1) {$v_3$};
\path (2.5,6) node     (y1) {$v_4$};\path (5.4,3.6) node     (y1) {$v_5$};

\end{tikzpicture}

\end{center}

  \caption{  A $(2,3)$-Cordial labellings of a 5-tournament}
\label{5tour}
\end{figure}
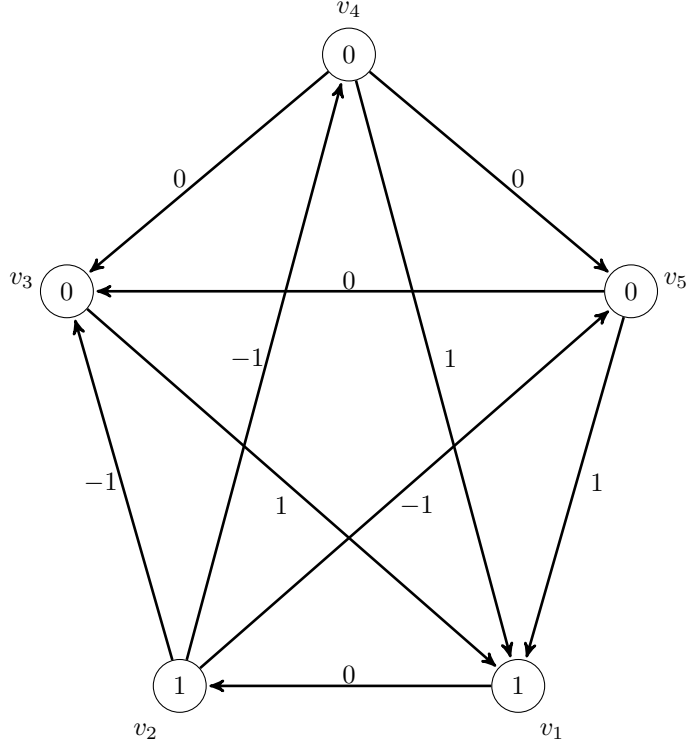

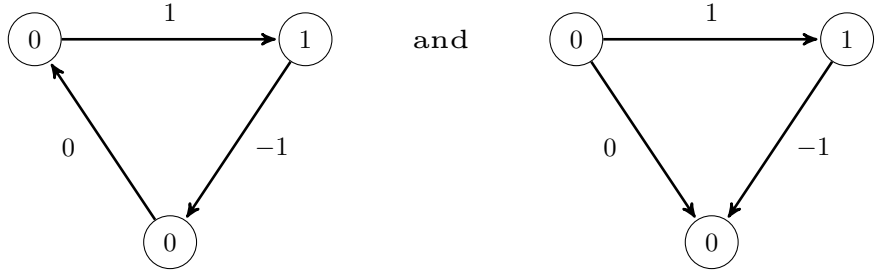
\begin{figure}
\begin{center}
\begin{tikzpicture}[scale=0.9]

\tikzset{vertex/.style = {shape=circle,draw,minimum size=2em}}
\tikzset{edge/.style = {->,> = stealth',shorten >=1pt,thick}}

\node[vertex] (v1) at  (0,6) {$0$};
\node[vertex] (v2) at  (2,3) {$0$};
\node[vertex] (v3) at  (4,6) {$1$};
\path (6,6) node  [scale=2]     (x) {
\tiny
 and} ;
\node[vertex] (v2k-2) at  (8,6) {$0$};
\node[vertex] (v2k-1) at  (10,3) {$0$};
\node[vertex] (v2k) at  (12,6) {$1$};

\draw[edge,->, line width=1.0pt] (v1) to (v3);
\draw[edge,->, line width=1.0pt] (v3) to (v2);
\draw[edge,->, line width=1.0pt] (v2) to (v1);
\draw[edge,->, line width=1.0pt] (v2k-2) to (v2k-1);
\draw[edge,->, line width=1.0pt] (v2k) to (v2k-1);
\draw[edge,->, line width=1.0pt] (v2k-2) to (v2k);

\path (2,6.4) node     (y1) {$1$};
\path (0.5,4.4) node     (y2) {$0$};
\path (3.5,4.4) node     (y3) {$-1$};
\path (8.5,4.4) node     (y4) {$0$};
\path (10,6.4) node     (y5) {$1$};
\path (11.5,4.4) node     (y6) {$-1$};

\end{tikzpicture}

\end{center}

  \caption{ $(2,3)$-Cordial labellings of 3- tournaments }
\label{3tour}
\end{figure}

\l  \label{5-tour} Every 5-tournament is $(2,3)$-cordial.\ll
\d  Let $T\in{\cal D}_5$ be a tournament.  Then there are two vertices, without loss of generality, $v_1$ and $v_2$, whose total out degree is four.  (And hence the total in-degree of  $v_1$ and $v_2$ is also four.  
  Let $f$ be the vertex labelling and let $f(v_1)=f(v_2)=1$ and $f(v_3)=f(v_4)=f(v_5)=0$.  Let $g$ be the arc labelling $g(\overrightarrow{v_iv_j})=f(v_j)-f(v_i)$.  Then, the arc between $v_1$ and $v_2$ is labelled 0, as are the three arcs between $v_3, v_4$ and $v_5$.  Thus there are four arcs labelled 0.  The three arcs from $v_1$ or $v_2$ to vertices $v_3, v_4$ or $v_5$  are labelled 1 and the three arcs from  $v_3, v_4$ or $v_5$ to vertices $v_1$ or $v_2$ are labelled -1.  In Figure \ref{5tour}   is an example of the labelling described above.   Thus $\Lambda_{f,g}({T})=(3,3,4)$.  That is $T$ is $(2,3)$-cordial.\dd

\l \label{6tour} If $n \geq 6$ and $T\in\Dn$ is a tournament on $n$ vertices then $T$ is not $(2,3)$-cordial.\ll
\d We divide the proof into two cases:

\noindent
{\bf Case 1.  $n$ is even.}  Let $n=2k$.  We shall show that there must be more arcs labelled 0 than is allowed in any $(2,3)$-cordial digraph with $\frac{n(n-1)}{2}$ arcs.  For any vertex labelled 0, there are $k-1$ other vertices also labelled 0 so that there are $k-1$ arcs labelled 0 that either begin or terminate at that vertex.  Also there are $k$ such vertices so there are $k(k-1)/2$ arcs between pairs of vertices each labelled 0.  (Note, since each arc is adjacent to two vertices we have divided the total number by 2 to get the number of distinct arcs labelled 0.)  There are also $k(k-1)/2$ arcs between pairs of vertices each labelled 1.  Thus we must have $k(k-1)$ arcs labelled 0.

Now, there must be at most one third the number of arcs labelled 0, so we must have $3k(k-1)\leq \frac{n(n-1)}{2}+2 = \frac{4k^2-2k+4}{2}$.  That is, we must have $k^2-2k-2 \leq 0$.  But that only happens if $k\leq 2$.  So if $k\geq 3$ or $n\geq 6$, $T$ is not $(2,3)$-cordial.

{\bf Case 2.  $n$ is odd.}  Let $n=2k+1$.   Without loss of generality, we may assume that there are $k$ vertices labelled 0 and $k+1$ vertices labelled 1.  Thus there are $\frac{1}{2} k(k-1)$ arcs labelled 0 that connect two vertices labelled  0 and $\frac{1}{2}(k+1)k$ arcs labelled 0 that connect two vertices labelled 1.   Thus there are at least $k^2$ arcs labelled 0.  To be $(2,3)$-cordial we must have that $3k^2\leq \frac{n(n-1)}{2}+2$, or $k^2-k-2\leq 0$.  That happens only if $k\leq 2$.  But, since $n$ is odd,  $n\geq 7$ so $k\geq 3$.  Thus, $T$ is not $(2,3)$-cordial. \dd

\begin{figure}
\begin{center}
\begin{tikzpicture}[scale=0.9]

\tikzset{vertex/.style = {shape=circle,draw,minimum size=2em}}
\tikzset{edge/.style = {->,> = stealth',shorten >=1pt,thick}}

\node[vertex] (v1) at  (0,6) {$ $};
\node[vertex] (v3) at  (4,6) {$ $};
\node[vertex] (v2) at  (0,3) {$ $};
\node[vertex] (v4) at  (4,3) {$ $};

\draw[edge,->, line width=1.0pt] (v3) to (v1);
\draw[edge,->, line width=1.0pt] (v3) to (v4);
\draw[edge,->, line width=1.0pt] (v4) to (v2);
\draw[edge,->, line width=1.0pt] (v2) to (v1);
\draw[edge,->, line width=1.0pt] (v2) to (v3);
\draw[edge,->, line width=1.0pt] (v4) to (v1);

\path (2.,2) node     (y33) {$T_{4,3}: (2,2,2,0)$};
\path (2,1.5)node (y44){Not (2,3)-cordial};

\node[vertex] (t1) at  (8,6) {$ $};
\node[vertex] (t3) at  (12,6) {$ $};
\node[vertex] (t2) at  (8,3) {$ $};
\node[vertex] (t4) at  (12,3) {$ $};

\draw[edge,->, line width=1.0pt] (t1) to (t3);
\draw[edge,->, line width=1.0pt] (t3) to (t4);
\draw[edge,->, line width=1.0pt] (t4) to (t2);
\draw[edge,->, line width=1.0pt] (t1) to (t2);
\draw[edge,->, line width=1.0pt] (t2) to (t3);
\draw[edge,->, line width=1.0pt] (t1) to (t4);

\path (10,2) node     (y33) {$T_{4,4}: (3,1,1,1)$};
\path (10,1.5)node (y44){Not (2,3)-cordial};

\node[vertex] (vv1) at  (0,12) {$0$};
\node[vertex] (vv3) at  (4,12) {$1$};
\node[vertex] (vv2) at  (0,9) {$1$};
\node[vertex] (vv4) at  (4,9) {$0$};

\draw[edge,->, line width=1.0pt] (vv1) to (vv3);
\draw[edge,->, line width=1.0pt] (vv3) to (vv4);
\draw[edge,->, line width=1.0pt] (vv4) to (vv2);
\draw[edge,->, line width=1.0pt] (vv2) to (vv1);
\draw[edge,->, line width=1.0pt] (vv3) to (vv2);
\draw[edge,->, line width=1.0pt] (vv1) to (vv4);

\path (2,12.4) node     (y1) {$1$};
\path (-.5,10.4) node     (y2) {$-1$};
\path (4.5,10.4) node     (y3) {$-1$};
\path (1.4,11.3) node     (y11) {$0$};
\path (2.6,11.3) node     (y22) {$0$};
\path (2,8.6) node     (y33) {$1$};

\path (2.,8) node     (y33) {$T_{4,1}: (2,2,1,1)$};

\node[vertex] (tt1) at  (8,12) {$0$};
\node[vertex] (tt3) at  (12,12) {$1$};
\node[vertex] (tt2) at  (8,9) {$0$};
\node[vertex] (tt4) at  (12,9) {$1$};

\draw[edge,->, line width=1.0pt] (tt1) to (tt3);
\draw[edge,->, line width=1.0pt] (tt3) to (tt4);
\draw[edge,->, line width=1.0pt] (tt4) to (tt2);
\draw[edge,->, line width=1.0pt] (tt1) to (tt2);
\draw[edge,->, line width=1.0pt] (tt3) to (tt2);
\draw[edge,->, line width=1.0pt] (tt1) to (tt4);

\path (10,12.4) node     (y1) {$1$};
\path (7.6,10.4) node     (y2) {$0$};
\path (12.5,10.4) node     (y3) {$0$};
\path (9.4,11.3) node     (y11) {$1$};
\path (10.54,11.3) node     (y22) {$-1$};
\path (10,8.6) node     (y33) {$-1$};

\path (10,8) node     (y33) {$T_{4,2}: (3,2,1,0)$};

\end{tikzpicture}

\end{center}

  \caption{ $(2,3)$-Cordial labellings of two 4-tournaments and two non  (2,3) cordial 4-tournaments with their out-degree sequences.}
\label{4-tours}
\end{figure}

\l\label{T4} The tournaments $T_{4,3}$ and $T_{4,4}$ of Figure \ref{4-tours} are not $(2,3)$-cordial.\ll\d Since $T_{4,4}$ is the reversal of $T_{4,3}$, by Lemma \ref{lab}  we only need show that $T_{4,3}$ is not  $(2,3)$-cordial.   Further, by Lemma \ref{lab} we may assume that the upper left vertex of $T_{4,3}$ in Figure \ref{4-tours} is labelled 0.  Since any permutation of the other three vertices results in an isomorphic graph we may assume that the upper right vertex is labelled 0 and the bottom two are labelled 1.  This results in one arc labelled 1, three arcs labelled -1 and two arcs labelled 0.  Thus $T_{4,3}$ is not  $(2,3)$-cordial.\dd

\t Let $T$ be an $n$-tournament.  Then $T$ is $(2,3)$-cordial if and only if $n\leq 5$ and $T$ is not isomorphic to $T_{4,3}$ or $T_{4,4}$.\tt
\d Lemmas \ref{T4},  \ref{5-tour} and \ref{6tour} together with Figures \ref{3tour} and \ref{4-tours}  establish the theorem.\dd

We end this section with a couple of observations we labell as corollaries:

\c The property of being (or not being) $(2,3)$-cordial is not closed under vertex deletion.\cc \d Every tournament on $k$ vertices is a vertex deletion of a tournament on $k+1$ vertices.  Thus, $T_{4,3}$, which is not  $(2,3)$-cordial, is a vertex deletion of a tournament on 5 vertices, which is  $(2,3)$-cordial, and this tournament is a vertex deletion of a tournament on 6 vertices, which is not  $(2,3)$-cordial.\dd

\c The property of being (or not being) $(2,3)$-cordial is not closed under arc contraction.\cc \d As in the above corollary, every tournament on $k$ vertices is an arc contraction  of a tournament on $k+1$ vertices.  Thus, $T_{4,3}$, which is not  $(2,3)$-cordial, is an arc contraction of a tournament on 5 vertices, which is  $(2,3)$-cordial, and this tournament is an arc contraction of a tournament on 6 vertices, which is not  $(2,3)$-cordial.\dd

%\pagebreak

\subsection{$(2,3)$-Orientations of Wheel Graphs}

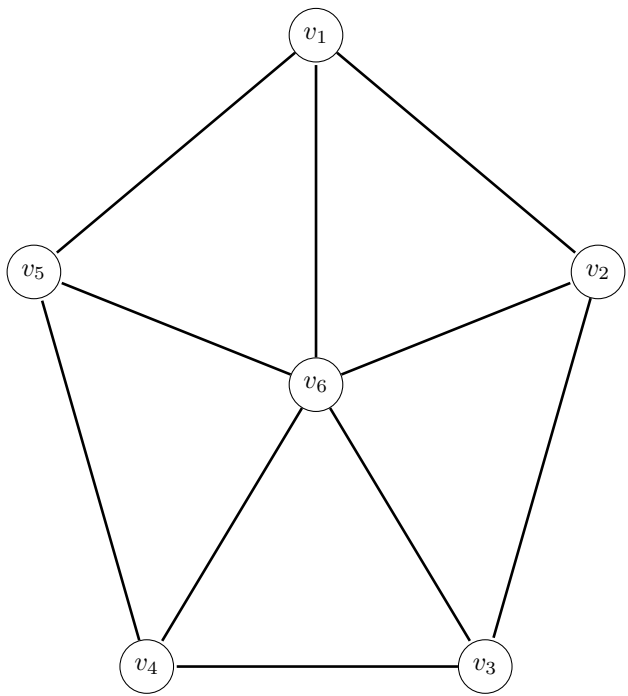
\begin{figure}[h]
\begin{center}
\begin{tikzpicture}[scale=1.5]

\tikzset{vertex/.style = {shape=circle,draw,minimum size=2em}}
\tikzset{edge/.style = {->,> = stealth',shorten >=1pt,thick}}

% vertices
\node[vertex] (v1) at  (2.5,5.6) {$v_1$};
\node[vertex] (v2) at  (5,3.5) {$v_2$};
\node[vertex] (v3) at  (4,0) {$v_3$};
\node[vertex] (v4) at  (1,0) {$v_4$};
\node[vertex] (v5) at  (0,3.5) {$v_5$};
\node[vertex] (v6) at  (2.5,2.5) {$v_6$};

\draw[edge,-, line width=1.0pt] (v1) to (v2);
\draw[edge,-, line width=1.0pt] (v2) to (v3);
\draw[edge,-, line width=1.0pt] (v3) to (v4);
\draw[edge,-, line width=1.0pt] (v4) to (v5);
\draw[edge,-, line width=1.0pt] (v1) to (v5);
\draw[edge,-, line width=1.0pt] (v6) to (v1);
\draw[edge,-, line width=1.0pt] (v6) to (v2);
\draw[edge,-, line width=1.0pt] (v6) to (v3);
\draw[edge,-, line width=1.0pt] (v6) to (v4);
\draw[edge,-, line width=1.0pt] (v6) to (v5);

\end{tikzpicture}

\end{center}

  \caption{  A 6-wheel graph.}
\label{6-wheel}
\end{figure}

A wheel graph on $n$ vertices consists of an $(n-1)$-star together with edges joining the non central vertices in a cycle.  A 6-wheel is shown in Figure \ref{6-wheel}.  Since we are not concerned with digraphs that contain digons, we shall assume that $n\geq 4$ in this section. 

\begin{figure}[h]
\begin{center}
\begin{tikzpicture}[scale=1.5]

\tikzset{vertex/.style = {shape=circle,draw,minimum size=2em}}
\tikzset{edge/.style = {->,> = stealth',shorten >=1pt,thick}}

% vertices
\node[vertex] (v1) at  (2.5,5.6) {$v_1$};
\node[vertex] (v2) at  (5,3.5) {$v_2$};
\node[vertex] (v3) at  (4,0) {$v_3$};
\node[vertex] (v4) at  (1,0) {$v_4$};
\node[vertex] (v5) at  (0,3.5) {$v_5$};
\node[vertex] (v6) at  (2.5,2.5) {$v_6$};

\draw[edge, ->, line width=1.0pt] (v1) to (v2);
\draw[edge,->, line width=1.0pt] (v2) to (v3);
\draw[edge,->, line width=1.0pt] (v3) to (v4);
\draw[edge,->, line width=1.0pt] (v4) to (v5);
\draw[edge,->, line width=1.0pt] (v5) to (v1);
\draw[edge,->, line width=1.0pt] (v6) to (v1);
\draw[edge,->, line width=1.0pt] (v6) to (v2);
\draw[edge,->, line width=1.0pt] (v6) to (v3);
\draw[edge,->, line width=1.0pt] (v6) to (v4);
\draw[edge,->, line width=1.0pt] (v6) to (v5);

\end{tikzpicture}

\end{center}

  \caption{  A 6-cycle-out-wheel graph.}
\label{6-wheel2}
\end{figure}
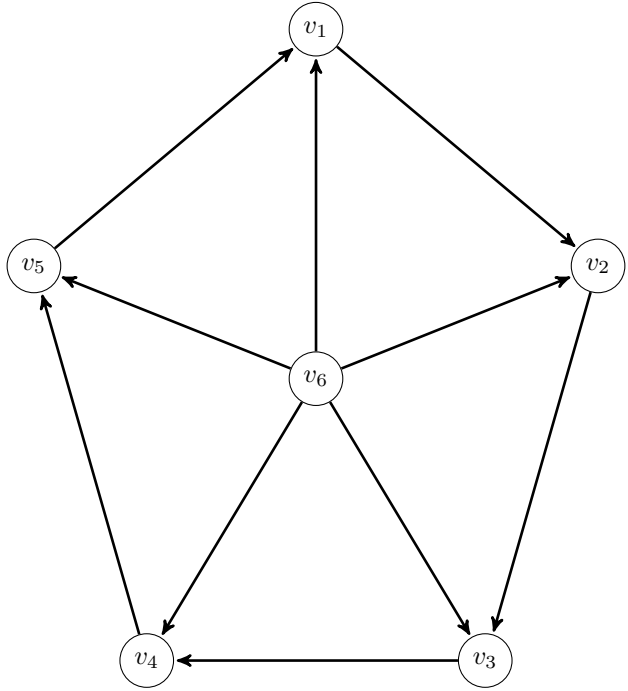

 An orientation of the wheel graph with the central vertex being a source/sink is called an out/in-wheel.  If the outer cycle of the wheel is oriented in a directed cycle the wheel is called a cycle-wheel.  If the $n$-wheel is oriented such that  it is both an out-wheel and a cycle-wheel it is called an $n$-cycle-out-wheel.  See Figure \ref{6-wheel2}.

\D Let $W_n=(V,A)$ be an $n$-wheel digraph with central vertex $v_n$ and with vertex labelling $f:V\to \{0,1\}$ Let $g$ be the induced arc labelling $g:A\to\{0,1,-1\}$ where $g(\overrightarrow{uv})=f(v)-f(u)$.  Let $S$ be the set of arcs incident with the central vertex and let $T$ be the set of arcs not incident with the central vertex.  
Define $\Lambda^S_{f,g}$ to be the real triple $\Lambda^S_{f,g}(D)=(\alpha_S,\beta_S,\gamma_S)$ where $\alpha_S=|g^{-1}(1)\cap S|, \beta_S=|g^{-1}(-1)\cap S|,$ and $\gamma_S =|g^{-1}(0)\cap S|$ and  $\Lambda^T_{f,g}$ to be the real triple $\Lambda^T_{f,g}(D)=(\alpha_{T},\beta_T,\gamma_T)$ where $\alpha_T=|g^{-1}(1)\cap T|, \beta_T=|g^{-1}(-1)\cap T|,$ and $\gamma_T =|g^{-1}(0)\cap T|$. Since $S\cup T= A$, the set of all arcs of $D$, $\alpha_S+\alpha_T=\alpha,  \beta_S+\beta_T=\beta$, and $\gamma_S+\gamma_T=\gamma$, where $\Lambda_{f,g}(D)=(\alpha,\beta,\gamma)$.\DD

\pagebreak

\t Let  $W_n$ be the $n$-wheel graph with central vertex $v_n$, and let $\overrightarrow W$ be a cyclic-out orientation of $W_n$.  Then $\overrightarrow W$ is not $(2,3)$-cordial.\tt
\d  We proceed with two cases, the case that $n$ is even then the case that $n$ is odd.   Let $W=(V,A)$ and let $f:V\to\{0,1\}$ be a vertex labelling and $g:A\to\{0,1,-1\}$ be the induced arc labelling, $g(\overrightarrow{uv}=f(v)-f(u)$.  Suppose that $f$ and $g$ is a $(2,3)$-cordial labelling. 

{\bf Case 1}:  $n=2k$.  Without loss of generality, we may assume that $f(v_n)=0$.  Thus, $\alpha_S=k,\beta_S=0$ and $\gamma_S=k-1$.  Further, since the orientation is cyclic, $\alpha_T=\beta_T$.  Since $\alpha-1\leq\beta\leq =\alpha+ 1$ we have $\alpha_T+k-1=\alpha_T+\alpha_S-1=\alpha-1\leq \beta = \beta_T+\beta_S= \beta_T=\alpha_T$.  Thus $k-1\leq 0$ or $k\leq1,$ a contradiction since $n\geq 4$.

{\bf Case 2}:  $n=2k+1$.  Without loss of generality we may assume that $f(v_n)=0$.  Since either $k$ or $k+1$ of the non central vertices must be labelled 0, we have two possibilities: $\alpha_S=k$ or $\alpha_S=k+1$.  

\noindent
{\bf Subcase 1}.  $\alpha_S=k$  Here we have $\gamma_S=k$ and $\beta_S=0$.  Further, $\alpha_T=\beta_T$ since $\overrightarrow W$ is cyclic.  Since the labelling is $(2,3)$-cordial, $\alpha-1\leq\beta\leq\alpha+1$.  Thus $k+\alpha_T-1=\alpha_S=\alpha_T-1= \alpha-1\leq\beta=\beta_S+\beta_T=\alpha_T$.  That is $k-1\leq0$ or $k\leq 1$, a contradiction since $n\geq4$.

\noindent
{\bf Subcase 2}.  $\alpha_S=k+1$  Here we have $\gamma_S=k-1$ and $\beta_S=0$.  Further, $\alpha_T=\beta_T$ since $\overrightarrow W$ is cyclic.  Since the labelling is $(2,3)$-cordial, $\alpha-1\leq\beta\leq\alpha+1$.  Thus $k+\alpha_T=k+1+\alpha_T-1=\alpha_S+\alpha_T-1= \alpha-1\leq\beta=\beta_S+\beta_T=\alpha_T$.  That is $k\leq0$, a contradiction since $n\geq4$.

In all cases we have arrived at a contradiction thus we must have that $\overrightarrow W$ is not $(2,3)$-cordial.\dd

\l  \label{rimarcs}Let $C$ be an undirected cycle with a $(0,1)$-vertex labelling.  Then, there is an even number of edges in $C$  whose incident vertices are labelled differently.\ll
\d  We may assume that the vertex $v_1$ is labelled 0.  Going around the cycle, the labelling goes from 0 to 1 then back again to zero.  This two step change must happen a fixed number of times then return to vertex $v_1$.  Thus there are an equal number of changes in labelling from 0 to 1 and from 1 to 0.  That is, the total number of changes is an even number.  \dd

\t \label{wn} Let $W_n$ be the undirected wheel graph on $n$ vertices.  Then, $W_n$ is not $(2,3)$-orientable if and only if $n=2k$ for some integer $k$, 4 does not divide $n$, and $2n-2=3z$ for some integer  $z$.  \tt

\d Let $\overrightarrow{W_n}$ be an orientation of the wheel graph on $n$ vertices with central vertex $v_n$.  Let $A_H$ be the set of arcs incident with $v_n$, and let $A_R$ be the arcs not incident with $v_n$.  Then $A=A_H\cup A_R$.  Let $f$ be a friendly vertex labelling and $g$ the induced arc labelling of $\overrightarrow{W_n}$.  Define , $\lambda_{f,H}(x)=|g^{-1}(x)\cap A_H|$, and $\lambda_{f,R}(x)=|g^{-1}(x)\cap A_R|$, Define $\lambda_f(x)=\lambda_{f,H}(x)+\lambda_{f,R}(x)$, that is $\lambda_{f}(x)=|g^{-1}(x)|$.  

We begin by showing that for  $n=2k$ for some integer $k$, $k=2\ell+1$ for some integer $\ell$, and $2n-2=3z$ for some integer  $z$ that $W_n$ is not $(2,3)$-orientable.  In this case, We may assume that $f(v_n)=0$ and $\lambda_{f,H}(1)+\lambda_{f,H}(-1)=k$, an odd integer.  By Lemma \ref{rimarcs} the number of arcs that are not incident with $v_n$ and  labelled either 1 or -1 is even.  Thus $\lambda_{f}(1)+\lambda_{f}(-1)=(\lambda_{f,H}(1)+\lambda_{f,H}(-1))+(\lambda_{f,R}(1)+\lambda_{f,R}(-1))$ is the sum of an even integer and an odd integer, so that  $\lambda_{f}(1)+\lambda_{f}(-1)$ is an odd integer.  But since the total number of arcs is $2n-2=3z$, if $\overrightarrow{W_n}$ is $(2,3)$-cordial, we must have $\lambda_{f}(1)+\lambda_{f}(-1)=2z$, an even integer,  a contradiction.  Thus, in this case $W_n$ is not $(2,3)$-orientable.

We now show that for all other cases $W_n$ is $(2,3)$-orientable.  We divide the proof into three cases, those being whether the total number of edges in $W_n$ is a multiple of three, one more than a multiple of three, or two  more than a multiple of three.

{\bf Case 1.} $2n-2=3z$ for some integer $z$.  

{\bf Subcase 1.1.} $n=2k$ and $k=2\ell$.    In this case, let $f$ be the labelling such that the labelling of the cycle has $2(z-\frac{k}{2})$ edges incident with vertices labelled differently.  Orient all arcs not incident with $v_n$ clockwise around the cycle.  Orient half the arcs incident with $v_n$ that are labelled 1 away from $v_n$,  and half toward $v_n$.  In this case, $\lambda_{f,H}(0)=k-1$, and $ \lambda_{f,H}(1)=\lambda_{f,H}(-1)=\frac{k}{2}=\ell$.  Further, $\lambda_{f,R}(1)=\lambda_{f,R}(-1)=z-\frac{k}{2}$.  Thus, $\lambda_{f}(1)=\lambda_{f}(-1)=\lambda_{f,H}(1)+\lambda_{f,R}(1)= \ell + z-\frac{k}{2}=z$.  Thus, we must also have $\lambda_f(0)=z$, and that $W_n$ is $(2,3)$-orientable.

{\bf Subcase 1.2.} $n=2k+1$.  In this case proceed as in Subcase 1 labelling the vertices not incident with $v_n$ with an even number of 1's (either $k$ or $k+1$).  Let $\ell$ be half of this even number.  Then, we produce a $(2,3)$-cordial orientation of $W_n$ the same as in Subcase 1.1.

{\bf Case 2.} $2n-2=3z+1$ for some integer $z$.  

{\bf Subcase 2.1.} $n=2k$, $k=2\ell$.   In this case, let $f$ be the labelling such that the labelling of the cycle has $2(z-\frac{k}{2})$ edges incident with vertices labelled differently.  Orient all arcs not incident with $v_n$ clockwise around the cycle.  Orient half the arcs incident with $v_n$ that are labelled 1 away from $v_n$,  and half toward $v_n$.   In this case, $\lambda_{f,H}(0)=k-1, \lambda_{f,H}(1)=\lambda_{f,H}(-1)=\frac{k}{2}=\ell$.  Further, $\lambda_{f,R}(1)=\lambda_{f,R}(-1)=z-\frac{k}{2}$.  Thus, $\lambda_{f}(1)=\lambda_{f}(-1)=\lambda_{f,H}(1)+\lambda_{f,R}(1)= \ell + z-\frac{k}{2}=z$.  Thus, we must also have $\lambda_f(0)=z+1$, and that $W_n$ is $(2,3)$-orientable.

{\bf Subcase 2.2.} $n=2k$, $k=2\ell+1$  In this case, let $f$ be the labelling such that the labelling of the cycle has $2(z-\frac{k-1}{2})$ edges incident with vertices labelled differently.    Orient $\ell$ of the arcs incident with $v_n$  that are labelled 1 away from $v_n$,  and $\ell+1$ of those arcs  toward $v_n$.    In this case, $\lambda_{f,H}(0)=k-1, \lambda_{f,H}(1)=\ell$ and $\lambda_{f,H}(-1)=\ell+1$.  Further, $\lambda_{f,R}(1)=\lambda_{f,R}(-1)=z-\frac{k}{2}=z=\ell$.  Thus, $\lambda_{f}(1)=\lambda_{f,H}(1)+\lambda_{f,R}(1)= \ell + z-\ell=z$, and  $\lambda_{f}(-1)=\lambda_{f,H}(-1)+\lambda_{f,R}(-1)= \ell+1 + z-\ell=z+1$.  Thus, we must also have $\lambda_f(0)= 2n-2-(z)-(z+1)=z$, and  thus  $W_n$ is $(2,3)$-orientable.

{\bf Subcase 2.3.} $n=2k+1$.  In this case proceed as in Subcase 2.1 labelling the vertices not incident with $v_n$ with an even number of 1's (either $k$ or $k+1$ depending upon whether $k$ is even or odd).  Let $\ell$ be half of this even number.  Then, we produce a $(2,3)$-cordial orientation of $W_n$ the same as in Subcase 2.1.

{\bf Case 3.} $2n-2=3z+2$ for some integer $z$.  

{\bf Subcase 3.1.} $n=2k$, $k=2\ell$.   In this case, let $f$ be the labelling such that the labelling of the cycle has $2(z-\frac{k}{2})$ edges incident with vertices labelled differently.  Orient all arcs not incident with $v_n$ clockwise around the cycle.  Orient half the arcs incident with $v_n$  that are labelled 1 away from $v_n$,  and half toward $v_n$.   In this case, $\lambda_{f,H}(0)=k-1, \lambda_{f,H}(1)=\lambda_{f,H}(-1)=\frac{k}{2}=\ell$.  Further, $\lambda_{f,R}(1)=\lambda_{f,R}(-1)=z-\frac{k}{2}$.  Thus, $\lambda_{f}(1)=\lambda_{f}(-1)=\lambda_{f,H}(1)+\lambda_{f,R}(1)= \ell + z-\frac{k}{2}=z$.  Thus, we must also have $\lambda_f(0)=z+1$, and that $W_n$ is $(2,3)$-orientable.

{\bf Subcase 3.2.} $n=2k$, $k=2\ell+1$  In this case, let $f$ be the labelling such that the labelling of the cycle has $2(z-\frac{k-1}{2})$ edges incident with vertices labelled differently.    Orient $\ell$ of the arcs incident with $v_n$  that are labelled 1 away from $v_n$,  and $\ell+1$ toward $v_n$.    In this case, $\lambda_{f,H}(0)=k-1, \lambda_{f,H}(1)=\ell$ and $\lambda_{f,H}(-1)=\ell+1$.  Further, $\lambda_{f,R}(1)=\lambda_{f,R}(-1)=z-\frac{k}{2}=z=\ell$.  Thus, $\lambda_{f}(1)=\lambda_{f,H}(1)+\lambda_{f,R}(1)= \ell + z-\ell=z$, and  $\lambda_{f}(-1)=\lambda_{f,H}(-1)+\lambda_{f,R}(-1)= \ell+1 + z-\ell=z+1$.  Thus, we must also have $\lambda_f(0)=z$, and that $W_n$ is $(2,3)$-orientable.

{\bf Subcase 3.3.} $n=2k+1$.  In this case proceed as in subcase 3.1 labelling the vertices not incident with $v_n$ with an even number of 1's (either $k$ or $k+1$ depending upon whether $k$ is even or odd.  Let $\ell$ be half of this even number.  Then, we produce a $(2,3)$-cordial orientation of $W_n$ the same as in Subcase 3.1.

We have now established the theorem.
\dd

\subsection{$(2,3)$-Orientations of Fan Graphs.}

A fan graph is isomorphic to a wheel graph with one edge of the cycle deleted.  Thus, by deleting one properly chosen arc from the cycle of a $(2,3)$-oriented $n$-wheel graph we obtain an orientation of the $n$-fan graph that is $(2,3)$-cordial.  Note that if there are at least as many  arcs labelled $x$ ($x=1,-1$ or $0$) as any other labelling, the properly chosen arc would be in the set of arcs labelled $x$.)  Thus there is only one case to consider, the case where $2n-2=3z, n=2k$ and $k=2\ell+1$ for some $z, k$, and $\ell$.

\begin{figure}[h]
\begin{center}
\begin{tikzpicture}[scale=1.26]

\tikzset{vertex/.style = {shape=circle,draw,minimum size=3em}}
\tikzset{edge/.style = {->,> = stealth',shorten >=1pt,thick}}

% vertices
\node[vertex] (v1) at  (0,3) {$v_1$};
\node[vertex] (v2) at  (0,6) {$v_2$};
\node[vertex] (v3) at  (1.5,6) {$v_3$};
\node[vertex] (v4) at  (3,6) {$v_4$};
\node[vertex] (v10) at  (6,6) {$v_{2\alpha-1}$};
\node[vertex] (v11) at  (7.5,6) {$v_{2\alpha}$};

\node[vertex] (v12) at  (7.5,4.5) {$v_{2\alpha+1}$};
\node[vertex] (v13) at  (7.5,1.5) {$v_{s-1}$};
\node[vertex] (v14) at  (7.5,0) {$v_{s}$};
\node[vertex] (v15) at (6,0) {$v_{s+1}$};
\node[vertex] (v16) at (1.5,0) {$v_{n-1}$};
\node[vertex] (v17) at (0,0) {$v_{n}$};

\path (4.5,6) node  (y1) {$\dots$};
\path (3.5,0) node  (y2) {$\cdots$};

\path (7.5,3) node  (z1) {$\vdots$};

\path (3.5,-.8) node  (y12) {$s=k+\alpha$};

\path (-.6,3) node  (z11) {0};
\path (1.5,6.6) node  (z11) {0};
\path (-.4,6.5) node  (z11) {1};
\path (3,6.6) node  (z11) {1};
\path (6,6.6) node  (z11) {0};
\path (7.8,6.5) node  (z11) {1};
\path (8.2,4.5) node  (z11) {1};
\path (8.2,1.5) node  (z11) {1};
\path (8.,-.34) node  (z11) {1};
\path (6,-.6) node  (z11) {0};
\path (1.5,-.6) node  (z11) {0};
\path (-.4,-.5) node  (z11) {0};

\draw[edge, ->, line width=1.0pt] (v1) to (v2);
\draw[edge,->, line width=1.0pt] (v2) to (v3);
\draw[edge,->, line width=1.0pt] (v3) to (v4);
\draw[edge,->, line width=1.0pt] (v4) to (y1);
\draw[edge,->, line width=1.0pt] (y1) to (v10);
\draw[edge,->, line width=1.0pt] (v10) to (v11);
\draw[edge,->, line width=1.0pt] (v11) to (v12);
\draw[edge,->, line width=1.0pt] (v12) to (z1);
\draw[edge,->, line width=1.0pt] (z1) to (v13);
\draw[edge,->, line width=1.0pt] (v13) to (v14);
\draw[edge,->, line width=1.0pt] (v14) to (v15);
\draw[edge,->, line width=1.0pt] (v15) to (y2);
\draw[edge,->, line width=1.0pt] (y2) to (v16);
\draw[edge,->, line width=1.0pt] (v16) to (v17);
\draw[edge,->, line width=1.0pt] (v17) to (v1);

\draw[edge, -, line width=.5pt] (v1) to (v3);
\draw[edge, -, line width=.5pt] (v1) to (v4);
\draw[edge, -, line width=.5pt] (v1) to (v10);
\draw[edge, -, line width=.5pt] (v1) to (v11);
\draw[edge, -, line width=.5pt] (v1) to (v12);
\draw[edge, -, line width=.5pt] (v1) to (v13);
\draw[edge, -, line width=.5pt] (v1) to (v14);
\draw[edge, -, line width=.5pt] (v1) to (v15);
\draw[edge, -, line width=.5pt] (v1) to (v16);

\end{tikzpicture}

\end{center}

 \caption{ $F_n$ with oriented cycle arcs.}
\label{nFan}
\end{figure}
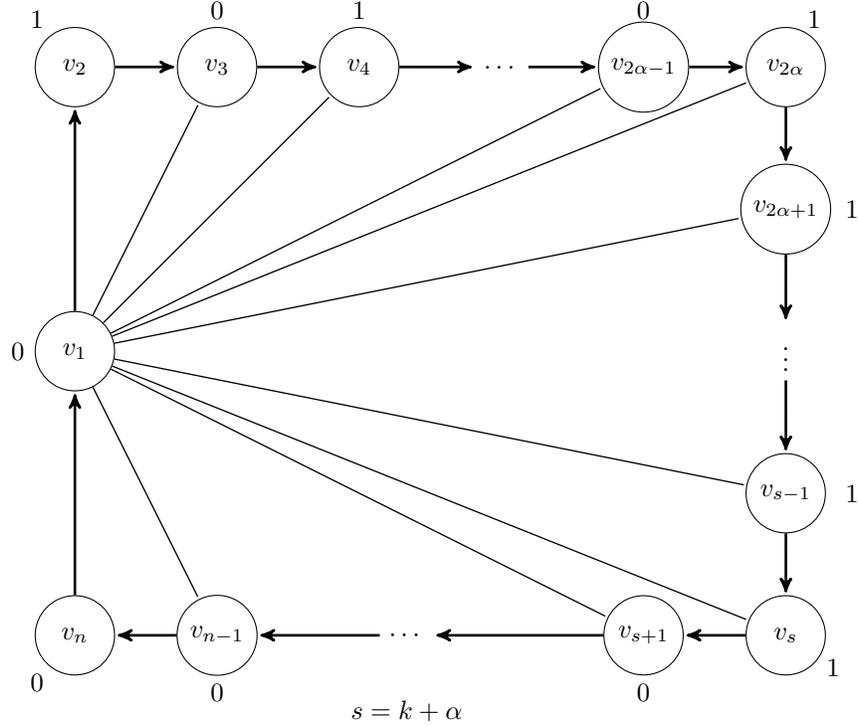

\t Let $n\geq 5$ and let  $F_n$ be the $n$-fan graph with central vertex $v_1$, that is the edges not on the cycle are all incident to $v_1$.  Let $\overrightarrow F$ be a cyclic-out orientation of $F_n$.  Then $\overrightarrow F$ is not $(2,3)$-cordial.\tt
\d As for wheel graphs, the number of arcs labelled 1 on the cycle is equal to the number of arcs labelled -1 and there are at least two arcs labelled 1 on the interior of the cycle.  thus, the number of arcs labled 1 in $\overrightarrow F$ is at least two more that the arcs labelled -1 in $\overrightarrow F$.  That is $\overrightarrow F$ is not $(2,3)$-cordial.\dd

\t  Let $F_n$ be the fan graph on $n$ vertices, $2n-3=3z+2$, $n=2k$ and $k=2\ell+1$ for some integers $k,\ell,$ and $z$.  Then there is an orientation of $F_n$ that is $(2,3)$-cordial.\tt
\d  Let $\alpha=z-\ell +1$, and define $f:V\to\{0,1\}$  by  $f(v_{2i-1})=0,$ $ i=1,\dots,\alpha$, $f(v_{2i})=1, i=1,\dots,\alpha$,   $f(v_{2\alpha+i})=1, i=1,\dots,k-\alpha$,   and $f(v_{k+\alpha+i})=0, i=1,\dots,k-\alpha$,   Note that $(k-\alpha)+(k+\alpha)=2k=n$, so all vertices are labelled.   Orient the cycle clockwise, so that the oriented cycle is $\overrightarrow{v_1v_2}$,   $\overrightarrow{v_2v_3}$, \dots,  $\overrightarrow{v_{n-1}v_n}$,  $\overrightarrow{v_nv_1}$.  See Figure \ref{nFan}  where the vertex labellings are outside the cycle.

Now, orient $\ell$ of the inner arcs from $v_1$ to arcs labelled 1 (except for $v_2$ which is not an inner arc) away from $v_1$ and the remaining $\ell$ such arcs inward so that we get $\overrightarrow{F_n}=(V,A)$.  Let $g:A\to\{0,1,-1\}$ be the induced labelling, $g(\overrightarrow{uv})=f(v)-f(u)$.  Then there are $\alpha$ arcs labelled $1$ on the cycle, $\alpha$ arcs labelled $-1$ on the cycle, $\ell$ of the inner arcs are labelled $-1$ and $\ell$ of the inner arcs are labelled $1$.  Thus, in all of  $\overrightarrow{F_n}$ there are $\alpha+\ell=z+1$ arcs labelled -1, $\alpha+\ell=z+1$ arcs labelled 1 and (hence) $z$ arcs labelled $0$.  That is, this orientation of $F_n$ is $(2,3)$-cordial.\dd

\subsection{Extremes of $(2,3)$ Cordiality}

As seen in section 3, complete graphs are not $(2,3)$-orientable if $n\geq 6$.  So the question arises:  How large  can a $(2,3)$-orientable graph be (how many edges)? Or: How large can a $(2,3)$-cordial digraph be?  That question was fully answered by M. A. Santana in \cite{SBMB}  For completeness we shall include the proofs of his results.

\t\cite[Theorem 3.1]{SBMB} \label{OR}Every simple directed graph is $(2,3)$ cordial if and only if there exists a friendly vertex labelling such that about $\frac{1}{3}$
of the edges are connected by vertices of the same label. 
\tt
\d %\begin{proof}
Let $G$ be a graph such that there exists a friendly labelling on $G$ such that about $\frac{1}{3}$ of the edges are connected by vertices of the same labell. This would mean about $\frac{2}{3}$ of the edges are connected by vertices of different labells, and therefore arcs may be assigned such that $G$ is cordial. Now let $H$ be a graph such that there does not exist a friendly labelling on $H$ such that such that about $\frac{1}{3}$ of the edges are connected by vertices of the same labell then there will be no way $H$ can be cordial since only then could about one third of the edges be labelled 0..
\dd %\end{proof}

Santana's application of Theorem \ref{OR} is 
\t\cite[Theorem 4.2]{SBMB} Given a directed graph $G = (V,E)$ with vertex set $V$ and $n = |V|$ with $n \geq 6$, and edge set $E$, the maximum size of $E$ such that $G$ is cordial for any given $n$ is
\tt
\begin{equation}
    |E|_{max} = \displaystyle{n\choose2} - Z + \Bigg \lceil \frac{1}{2}\Big ( \displaystyle{n\choose2} - Z\Big)  \Bigg \rceil\\
    Z \\= \displaystyle{\lceil \frac{n}{2}  \rceil\choose 2} + \displaystyle{\lfloor \frac{n}{2}  \rfloor\choose 2}.
\end{equation}
\begin{figure}
    \centering
    \begin{tikzpicture}[scale=0.65]
     \node[shape=circle,draw=black] (A) at (2,-4) {0};
     \node[shape=circle,draw=black] (B) at (-3,-4) {0};
     \node[shape=circle,draw=black] (C) at (-5,0) {0};
     \node[shape=circle,draw=black] (D) at (-2,4) {1};
     \node[shape=circle,draw=black] (E) at (2,4) {1};
     \node[shape=circle,draw=black] (F) at (5,0) {1};
     \draw[dashed] (A) -- (B);
     \draw[dashed] (A) -- (C);
     \path [-] (A) edge node[left] {} (D);
     \path [-] (A) edge node[left] {} (E);
     \path [-] (A) edge node[left] {} (F);
     
     \draw[dashed] (B) -- (C);
     \path [-] (B) edge node[left] {} (D);
     \path [-] (B) edge node[left] {} (E);
     \path [-] (B) edge node[left] {} (F);
     
     \path [-] (C) edge node[left] {} (D);
     \path [-] (C) edge node[left] {} (E);
     \path [-] (C) edge node[left] {} (F);
     
     \draw[dashed] (D) -- (E);
     \draw[dashed] (D) -- (F);
     
     \draw[dashed] (F) -- (E);
\end{tikzpicture}
\caption{A complete graph. Dashed Lines represent edges labelled zero regardless of arc orientation}
\end{figure}
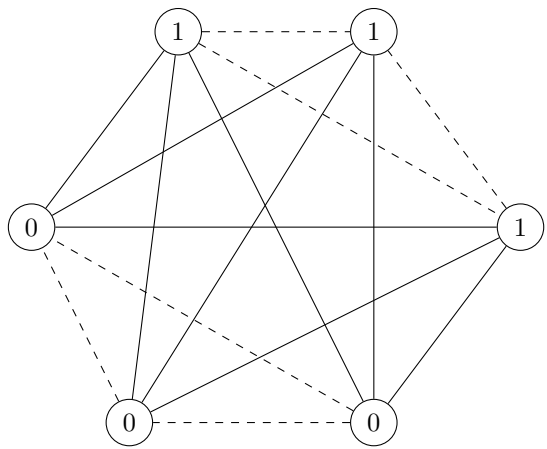
\d %\begin{proof}
%\todo{cite the other paper}
From section 3 we have that for any tournament with $n \leq 5$ there exists a cordial labelling, save for the case when $n =4$ thus we begin with a complete graph with $n \geq 6$. Recall that the number of vertices on a complete graph is $n\choose2$. Due to our cordial labelling the number of edges with an induced labelling of 0 will be our $Z$. This is because it will be the number of edges connected by two vertices of the same label, as shown in Figure 7. If $n$ is even that will mean that $Z = 2 {\frac{n}{2}\choose2}$, i.e., it will be the number of edges on two $\frac{n}{2}$ complete graphs represented by the labellings of ones and zeros. The floor and ceiling function in $(1)$ simply account for the odd case.\\
For every tournament with $n \geq 6$ vertices, $Z > \frac{1}{3} {n\choose 2}$. Therefore some of the arcs labelled zero will need to be removed to get a cordial graph. How many arcs need to be removed is going to be equal to how much greater $Z$ is than the number of half the number of arcs not labelled zero. By the definition of a directed cordial graph we know that $Z$ can be larger than $\alpha$ or $\beta$ and we can still have a cordial graph, hence the ceiling function.
\dd %\end{proof}

As mentioned in the introduction, the smallest non $(2,3)$-cordial digraph is an orientation of {  X}$_n$, three parallel arcs.  A question may be asked:  What is the minimum number of arcs in a non $(2,3)$-cordial digraph that has no isolated vertices?

\section{Conclusions}
In this article we have shown that the only tournaments that are $(2,3)$-cordial are when $n\leq 5$ and then not for two  4-tournaments.    That except for one case when $n$ is even, the $n$-wheel graph has an orientation that is $(2,3)$-cordial and that at least one orientation of any wheel graph is not $(2,3)$-cordial.  Further, we show that every fan graph has a $(2,3)$-cordial orientation, and as for wheel graphs there is always an orientation of the $n$-fan that is not $(2,3)$-cordial.

%\pagebreak

Authors email addresses: {LeRoy B. Beasley,  leroy.b.beasley@aggiemail.usu.edu; Manuel Santana, manuelarturosantana@gmail.com;    Jonathan Mousley,  jonathanmousley@gmail.com; David E. Brown,  david.e.brown@usu.edu}

\end{document}